\documentclass[journal]{IEEEtran}

\usepackage{cite}

\usepackage{geometry}
\usepackage{cite}

\ifCLASSINFOpdf
\else
\fi
\usepackage{array}

\usepackage{amssymb,amsmath}
\usepackage{graphicx}
\usepackage{amsthm}
\usepackage{epstopdf}
\allowdisplaybreaks[4]
\usepackage{flushend}
\usepackage{color}
\usepackage{fancyhdr}

\emergencystretch=\maxdimen
\newtheorem{theorem}{Theorem}
\newtheorem{definition}{Definition}

\newtheorem{lemma}{Lemma}

\newtheorem{assumption}{Assumption}
\newtheorem{remark}{Remark}
\newtheorem{proposition}{Proposition}

\newcommand{\col}{$\upshape{col}$}
\newcommand{\blk}{$\upshape{blk}$}

\begin{document}
%
\title{Distributed continuous-time strategy-updating rules for noncooperative games with discrete-time communication}
%
%
%

\author{\IEEEauthorblockN{Xin~Cai, Feng~Xiao, {\it Member, IEEE}, Bo~Wei and Fang~Fang, {\it Senior Member, IEEE}}

\thanks{This work has been submitted to the IEEE for possible publication. Copyright
may be transferred without notice, after which this version may no longer be
accessible. This work was supported in part by the National Natural Science Foundation of China (NSFC, Grant No. 62273145), in part by Beijing Natural Science Foundation (Grant No.~4222053), and in part by the Natural Science Foundation of Xinjiang Uygur Autonomous Region (Grant No. 2022D01C694). (\textit{Corresponding author: Feng~Xiao})

X. Cai and F. Xiao are with the State Key Laboratory of Alternate Electrical Power System with Renewable Energy Sources and with the School of Control and Computer Engineering, North China Electric Power University, Beijing 102206, China (emails: caixin\_xd@126.com; fengxiao@ncepu.edu.cn). X. Cai is also with the School of Electrical Engineering, Xinjiang University, Urumqi 830047, China.

B. Wei and F. Fang are with the School of Control and Computer Engineering,
North China Electric Power University,
Beijing 102206, China (emails: bowei@ncepu.edu.cn; ffang@ncepu.edu.cn).}}

%
%

\markboth{This work has been submitted to the IEEE for possible publication.} 
{}

%



\maketitle

\begin{abstract}
In this paper, continuous-time noncooperative games in networks of double-integrator agents are explored. The existing methods require that agents  communicate with their neighbors in real time. In this paper, we propose two discrete-time communication schemes based on the designed continuous-time strategy-updating rule for the efficient use of communication resources. First, the property of the designed continuous-time rule is analyzed to ensure that all agents' strategies can converge to the Nash equilibrium. Then, we propose periodic and event-triggered communication schemes for the implementation of the designed rule with discrete-time communication. The rule in the periodic case is implemented synchronously and easily. The rule in the event-triggered case is executed asynchronously without Zeno behaviors. All agents in both cases can asymptotically reach to the Nash equilibrium by interacting with neighbors at discrete times. Simulations are performed in networks of Cournot competition to illustrate the effectiveness of the proposed methods.
\end{abstract}


%
\IEEEpeerreviewmaketitle

\section{Introduction}
In recent years, there has been a growing interest in adopting game theory to characterize the interactions among decision-making agents in distributed control systems \cite{Ratliff.2016}. By this method, the coordination objective of multi-agent systems can be formulated as the Nash equilibrium (NE) which corresponds to the natural emergent behavior resulting from the interactions of selfish players \cite{Cortes.2015}.

Roughly speaking, the dynamics of noncooperative games can be classified into two categories. One relies on the best-response dynamics, which derive from the most natural ``game playing", i.e., each player changes his strategy to minimize his cost, given strategies of other players \cite{Nisan.2007, Cortes.2015, Swenson.2018,Chen.2019}. The other relies on the gradient dynamics, which are usually suitable for continuous cost functions. From the perspective of optimization, the strategy of each player is updated along the direction of maximal descent of his cost function.

Based on the gradient dynamics, the studies of continuous-time noncooperative games in the framework of multi-agent systems have recently attracted an increasing attention. In this framework, the distributed strategy-updating rules (or NE seeking algorithms) based on consensus protocols designed for games with incomplete information are totally different from the methods proposed in Bayesian games\cite{Zheng.2018,Qu.2018,Yu.2018}. In such a setting, agents are forced to move to NE with the estimation of other agents' strategies by communicating with neighbors. The corresponding results are more suitable for networked control systems, such as sensor networks \cite{Stankovic.2012, Ma.2019}, communication networks \cite{Gharesifard.2013}, mechanical systems \cite{Deng.2019} and smart grids \cite{Ma.2013}.  In consideration of the cyber-physical scenarios where the implementation of distributed algorithms among agents is driven both by inherent dynamics and their own interests, some studies have focused on the games among agents with complex dynamics. A distributed NE seeking algorithm was proposed for mobile robots with second-order dynamics in the sensor coverage\cite{Ibrahim.2018}. A network of Euler-Lagrange systems was steered by distributed algorithms to NE of aggregative games\cite{Deng.2019}. Besides, NE seeking algorithms were proposed for agents with disturbed multi-integrator dynamics \cite{Romano.2020}. Multi-integrator agents were also studied in noncooperative games with coupled constraints\cite{Bianchi.2019}. A network of passive nonlinear second-order systems was steered by a distributed feedback algorithm to the Cournot-Nash equilibrium\cite{Persis.2019c}. Zhang \textit{et al.} described players in aggregative games as nonlinear dynamic systems perturbed by external disturbances and designed a distributed algorithm for agents to arrive at the NE and to reject disturbances simultaneously\cite{Zhang.2019}. Moreover, the seeking algorithm was proposed for quadratic games with nonlinear dynamics \cite{Huang.2020}.

However, a common feature of continuous-time seeking algorithms in all aforementioned literature is that the communications among agents are continuous. For example, in mobile sensor networks, agents (e.g., robots or vehicles) are forced by continuous-time seeking algorithms\cite{Bianchi.2019, Romano.2020}, which require to be implemented by continuous communications among agents. Taken the communication costs into account, the continuous-time mechanism is undoubtedly expensive. Furthermore, due to the limited storage energy in mobile sensors, it is impossible to apply these algorithms in practical situations where agents communicate with each other only at discrete-time instants. To the best of our knowledge,  only a few work has been reported on discrete-time communication for noncooperative games. Furthermore, few studies have been reported on the discrete-time interactions in continuous-time NE seeking algorithms. A self-triggered communication scheme was designed in best-response dynamics \cite{Cortes.2015}. Static event-triggering schemes were designed for differential games \cite{Yuan.2018,Xue.2020}. These observations motivate us to study the strategy-updating rules with discrete-time communication for noncooperative games. The interaction laws among agents are designed to accurately estimate the strategies of others. During the time-interval between any two successive discrete communication times, the outdated information is used. The longer the communications are disconnected, the further away the estimations could be driven from the real values. It may happen that agents' strategies are far away from the NE of games. In such a case, how to design discrete-time communication schemes to ensure the convergence of all agents' strategies to the NE is a critical problem to be solve. Thanks to the extensively studied the control and communication technology based on event-triggered schemes  \cite{Girard.2015, Li.2016,Liu.2019,Sun.2020,Ge.2020}, they make it possible to solve this problem.

This paper is to investigate discrete-time communication schemes of a distributed continuous-time strategy-updating rule for agents who have double-integrator dynamics. Similar to the existing work, all agents need to agree to implement the designed algorithm by exchanging information with their neighbors. Thus, agents can estimate the strategies of others through local interactions in a distributed way. The main contributions of this paper are summarized as follows.

1) To remove the requirement of the continuous communication in the existing distributed NE seeking algorithms, a periodic communication scheme is firstly proposed to realize periodic and synchronous sampling. Then, for the purpose of saving communication resources, we design a asynchronously dynamic event-triggered communication scheme which depends on an internal variable with its own dynamics. The scheme is proofed to be Zeno-free. And it is different from the self-triggered communication \cite{Cortes.2015} and the static event-triggering scheme\cite{Yuan.2018,Xue.2020}. The self-triggered scheme schedules the next event time to sample new information by the current triggering time and available information. But the event-triggering scheme broadcasts new information at event-triggering time only determined by local information.

2) Agents with double-integrator dynamics are studied in noncooperative games. Different from the updating rule for the predicted strategy containing positions and velocities proposed in\cite{Romano.2020}, the strategy-updating rule given in this paper utilizes the agent's own real time strategy (i.e., positions) and the estimation of strategies of other agents for forcing itself to update its strategy along the direction of the subgradient of its cost function.

3) Continuous cost functions are investigated in this paper, and the assumption on continuous differentiable cost functions in related previous works is relaxed \cite{Ye.2017, Deng.2019, Gadjov.2019, Lu.2019, Cai.2020}. The proposed strategy-updating rules synthesize subgradient dynamics and differential inclusions to deal with continuous cost functions.

The rest of this paper is organized as follows. In Section \uppercase\expandafter{\romannumeral2}, the problem formulation is given and some related preliminaries are introduced. In Section \uppercase\expandafter{\romannumeral3}, a continuous-time strategy-updating rule is proposed. The discrete-time communication schemes are designed in Section \uppercase\expandafter{\romannumeral4}. Simulation examples are provided in Section \uppercase\expandafter{\romannumeral5}. Finally, some conclusions and future topics are stated in Section \uppercase\expandafter{\romannumeral6}.

Notations: $\mathbb{R}$ denotes the set of real numbers. $\mathbb{R}_{\geq0}$ is the set of non-negative real numbers. $\mathbb{Z}_{\geq0}$ is the set of non-negative integer numbers. $\mathbb{R}^n$ is the $n$-dimensional real vector space. $\mathbb{R}^{n\times m}$ denotes the set of $n \times m$ real matrices. Given a vector $x\in\mathbb{R}^n$, $\|x\|$ is the Euclidean norm. $A^T$ and $\|A\|$ are the transpose and the spectral norm of matrix $A\in \mathbb{R}^{n\times n}$, respectively. For matrices $A$ and $B$, $A\otimes B$ denotes their Kronecker product. $\lambda_{2}(A)$ and $\lambda_n(A)$ are the second smallest and the largest eigenvalues of matrix $A$, respectively, and they are expressed simply as $\lambda_{2}$ and $\lambda_n$ . Let $\col(x_1,\ldots,x_n)=[x_1^T,\ldots,x_n^T]^T$.  $\blk\{A_1,\ldots,A_n\}$ is a block diagonal matrix with diagonal elements $A_1, \ldots, A_n$. $\boldsymbol{1}_n$ and $\boldsymbol{0}_n$ are $n$-dimensional column vectors where all elements are 0 and 1, respectively. $I_n$ denotes the $n\times n$ identity matrix. A zero matrix is denoted by $\boldsymbol{0}$ with an appropriate dimension.  A set-valued map $\mathcal{F}(x):\mathbb{R}^n\rightrightarrows \mathbb{R}^n$ is the map from a vector $x\in\mathbb{R}^n$ to the collection of all subsets of $\mathbb{R}^n$.

\section{Problem formulation and Preliminaries}

\subsection{Problem Formulation}
An $N$-person noncooperative game with $N$ agents is considered here. Let $G$ $=$ $(\mathcal{I},\Omega,J)$ denote the game, where agents are indexed in the set $\mathcal{I}$ $=$ $\{1,\ldots,N\}$,  $\Omega$ $=$ $\Omega_1\times\cdots\times\Omega_N\subset\mathbb{R}^{Nn}$ is the strategy space of the game with the strategy set $\Omega_i\subset \mathbb{R}^{n}$ of agent $i\in \mathcal{I}$. $J$ $=$ $(J_1,\ldots,J_N)$, where $J_i(x_i,x_{-i})$ $:$ $\Omega_i\times \prod_{j\neq i}\Omega_j\rightarrow \mathbb{R}$ is agent $i$'s cost function depending on its own strategy $x_i$ $\in$ $\Omega_i$ and the other agents' strategies denoted by a vector $x_{-i}$ $=$ $\col(x_1,\ldots,x_{i-1},x_{i+1},\ldots,x_N)$. Let $x$ $=$ $\col(x_1,\ldots,x_N)\in \Omega$, which denotes the strategy profile composed of all agents' strategies. For the game with incomplete information, agents have to estimate strategies of the others by interactions with their neighbors in a network $\mathcal{G}$ $=$ $(\mathcal{I},\mathcal{E})$, which may be undirected or directed. For the detailed concepts of graphs, please refer to \cite{Godsil.2001}.

We consider that each agent in the network has the inherent double-integrator dynamics,
\begin{equation} \label{ad}
\left\{\begin{array}{l}
\dot{x}_i=v_i,\\
\dot{v}_i=u_i,
\end{array}\right.
\end{equation}
where $x_i\in\mathbb{R}^{n}$, $v_i\in\mathbb{R}^n$, and $u_i\in\mathbb{R}^{n}$ are the agent $i$'s strategy, auxiliary state, and control input, respectively. Given the other agents' strategies $x_{-i}$, agent $i$ aims to minimize its cost in the game, i.e., $\min_{x_i\in \Omega_i}J_i(x_i,x_{-i})$. The double-integrator dynamics \eqref{ad} can represent several types of physical systems, such as mobile robots in sensor networks\cite{Stankovic.2012}, autonomous vehicles in traffic scenarios \cite{Dreves.2018} and Euler-Lagrange systems\cite{Deng.2019}.

\begin{assumption} \label{as1}
Let $\Omega_i=\mathbb{R}^n$. The cost function $J_i(x_i,x_{-i})$ is continuous in all its arguments $x$ and convex in $x_i$ for every fixed $x_{-i}$ and for all $i\in\mathcal{I}$.
\end{assumption}
Under Assumption 1, the $N$-person noncooperative game formulated in this paper admits a NE \cite[Theorem 4.4]{Basar.1999}. The cost function $J_i$ represents a goal or performance metric of agent $i$, and it may be nonsmooth in many settings. For example, a piecewise linear price function was studied in the Cournot model \cite{Hobbs.2007}, the performance of compressing sensing is measured by $l_1$-norm in \cite{Jafarpour.2011}, and the congestion costs of flow control are assumed to be piecewise smooth in communication networks in \cite{Yin.2011}.

In summary, the problems we need to handle include: 1) the proposal of continuous-time strategy-updating rule for double-integrator agents with continuous cost functions; 2) the design of discrete-time communication schemes to reduce communication loads; 3) the analysis of the proposed methods that can ensure the asymptotical convergence of strategies to the unique NE of noncooperative games.

\begin{remark}
Assumption \ref{as1} is more general than that in \cite{Cai.2020}. Under Assumption \ref{as1}, it is seen that the strategy-updating rule designed in \cite{Cai.2020} is a special case corresponding to the rule designed in this paper. Moreover, the analysis in \cite{Cai.2020} is only based on ordinary differential equations and classical Lyapunov stability theory. There is no further study on the periodic communication scheme in the conference version of this paper \cite{Cai.2020}. Although a part of conclusions in this paper are similar to those in \cite{Cai.2020}, the conclusions in this paper are applicable to the cases with more general cost functions.
\end{remark}

\subsection{Preliminaries}
Here, we introduce some necessary notations and lemmas in noncooperative games, convex analysis and differential inclusions.

\begin{definition}[\!\!{\cite[Definition 3.7]{Basar.1999}}]
A pure NE of game $G=(\mathcal{I},\Omega,J)$ is a strategy profile $x^*=\col(x_1^*,\ldots,x_{N}^*)\in \Omega$ satisfying the following inequality
\begin{align*}
J_i(x_i^*,x_{-i}^*)\leq J_i(x_i,x_{-i}^*)
\end{align*}
for any $x_i\in \Omega_i$ and all $i\in\mathcal{I}$.
\end{definition}
The NE is the point where any agent has no willingness to decrease its cost by changing its strategy unilaterally..

\begin{lemma}[\!\!{\cite[Corollary 4.2]{Basar.1999}}] \label{lemma1}
Under Assumption \ref{as1}, the game $G=(\mathcal{I},\Omega,J)$ admits a pure NE $x^*\in \Omega$ satisfying
\begin{equation} \label{oc}
\boldsymbol{0}_n\in\partial_{x_i} J_i(x_i^*,x_{-i}^*),  \ \forall i\in \mathcal{I},
\end{equation}
where $\partial_{x_i} J_i(x_i,x_{-i})\in\mathbb{R}^n$ is the subdifferential of cost function $J_i$ at the strategy $x_i$ for fixed $x_{-i}$.
\end{lemma}

The condition \eqref{oc} is a generalization of $\nabla_iJ_i(x_i^*,x_{-i}^*)=\boldsymbol{0}_n$ which denotes the gradient of $J_i(x_i,x_{-i})$, which is continuous differential, with respect to $x_i$ and is widely used in \cite{Basar.1999, Gharesifard.2013, Ye.2017, Deng.2019, Romano.2020, Gadjov.2019}. Let $F(x)=\prod_{i=1}^N\partial_{x_i}J_i(x_i,x_{-i})$.

\begin{assumption} \label{as4}
The map $F(x):\Omega\rightrightarrows \mathbb{R}^{Nn}$ is Lipschitz continuous with Lipschitz constant $\theta>0$ and is strongly monotone; that is, $(x-x')^T(d-d')>w\|x-x'\|^2$, $w>0$, $\forall x,x'\in\Omega$, $d\in F(x)$, $d'\in F(x')$.
\end{assumption}

\begin{definition}[\!\!{\cite[Theorem 6.1.2]{Jean.2001}}]
If a function $f(\cdot): C \rightarrow \mathbb{R}$ is strongly convex on a convex set $C\subset \mathbb{R}^n$, there exists a constant $c>0$ such that
\begin{align*}
(y-y')^T(d-d')\geq c\|y-y'\|^2, \ \forall d\in\partial f(y), d'\in\partial f(y'),
\end{align*}
where $\partial f(y)$ is the subdifferential of $f(\cdot)$ at $y$.
\end{definition}

 \begin{proposition}[\!\!{\cite[Proposition 9]{Cortes.2008}}]
 Let $f(y):\mathbb{R}^n \rightarrow \mathbb{R}$ be a locally Lipschitz and convex function. Then,

 (1) $\partial f(y)\subset \mathbb{R}^n$ is nonempty, convex and compact, and all $d\in \partial f(y) $ satisfy that $\|d\|\leq l$ for some $l>0$;

 (2) $\partial f(y)$ is upper semi-continuous at $y\in \mathbb{R}^n$.
 \end{proposition}

A differential inclusion \cite{Aubin.1984} is considered as follows
\begin{equation} \label{di}
\dot{z}\in \mathcal{F}(z), z(0)=z_0,
\end{equation}
where $\mathcal{F}:\mathbb{R}^n\rightrightarrows \mathbb{R}^n$ is a set-valued map. A solution of \eqref{di} is an absolutely continuous curve $z:[0,T]\rightarrow \mathbb{R}^n$ that satisfies \eqref{di} for almost all $t\in [0,T]$. The set of equilibria of the system \eqref{di} is $E=\{z\in \mathbb{R}^n| 0\in \mathcal{F}(z)\}$.
If $\mathcal{F}: \mathbb{R}^n\rightrightarrows \mathbb{R}^n$ is upper semicontinuous with nonempty, compact, and convex values, there exists a solution to \eqref{di} for any initial condition\cite[Lemma 2.3]{Gharesifard.2013}.

\begin{lemma}[\!\!{\cite[Theorem 4]{Cortes.2008}}] \label{lemma2}
Let $V:\mathbb{R}^n \rightarrow \mathbb{R}$ be a smooth function. Define $S\subset\mathbb{R}^n$ as a strongly positively invariant set under \eqref{di}. The set-valued Lie derivative of $V$ with respect to $\mathcal{F}$ at $z$ is
\begin{align*}
\mathcal{L}_{\mathcal{F}}V=\{\zeta\in\mathbb{R}|\zeta=(\nabla V(z))^T \nu, \forall \nu \in \mathcal{F}(z)\}.
\end{align*}
If $max \mathcal{L}_{\mathcal{F}}V\leq0$ or $\mathcal{L}_{\mathcal{F}}V=\emptyset$, $\forall z\in S$, and the evolutions of \eqref{di} are bounded, the solutions of \eqref{di} starting from $S$ converge to the largest weakly positively invariant set $\mathcal{M}$ contained in $S\cap \{z\in\mathbb{R}^n|0\in\mathcal{L}_{\mathcal{F}}V\}$. When $\mathcal{M}$ is finite, the limit of every solution exists and is an element of $\mathcal{M}$.
\end{lemma}

\section{Distributed continuous-time strategy-updating rule}
For the games with incomplete information, every agent is assumed to estimate the other agents' strategies and to regulate the estimation by communication with its neighbors.

Inspired by augmented pseudo-gradient dynamics designed in \cite{Gadjov.2019}, we define $\boldsymbol{x}^i=\col(x_1^i,\ldots,x_N^i)$ as the estimation vector of agent $i$ about all agents' strategies, where $x_j^i$ is the player $i$'s estimation on the strategy of player $j$, and $x_i^i=x_i$. Denote $\boldsymbol{x}_{-i}^i=\col(x_1^i,\ldots,x_{i-1}^i,x_{i+1}^i,\ldots,x_N^i)$ and let $\alpha$ be a positive constant to be designed. With the estimation vector $\boldsymbol{x}_{-i}^i$, $\partial_{x_i}J_i(x_i,\boldsymbol{x}_{-i}^i)$ represents the subdifferential  of $J_i(x_i,\boldsymbol{x}_{-i}^i)$ at $x_i$. We propose the following strategy-updating rule.
\begin{equation} \label{csur}
\begin{aligned}
\dot{x}_i&=v_i,\\
\dot{v}_i&\in-kv_i-\partial_{x_i}J_i(x_i,\boldsymbol{x}_{-i}^i)-\frac{\alpha}{k}R_i\sum_{j\in\mathcal{I}}a_{ij}(\boldsymbol{x}^i-\boldsymbol{x}^j),\\
\dot{\boldsymbol{x}}_{-i}^i&=-\alpha S_i\sum_{j\in \mathcal{I}}a_{ij}(\boldsymbol{x}^i-\boldsymbol{x}^j),
\end{aligned}
\end{equation}
where $a_{ij}$ is the weight on edge $(i,j)$ $\in$ $\mathcal{E}$ of graph $\mathcal{G}$.
\begin{align*}
R_i&=[\boldsymbol{0}_{n\times(i-1)n}, I_n, \boldsymbol{0}_{n\times(N-i)n}]
\end{align*}
and
\begin{align*}
S_i&=\begin{bmatrix}I_{(i-1)n} & \boldsymbol{0}_{(i-1)n\times n} & \boldsymbol{0}_{(i-1)n\times(N-i)n}\\ \boldsymbol{0}_{(N-i)n\times(i-1)n} & \boldsymbol{0}_{(N-i)n\times n} & I_{(N-i)n} \end{bmatrix}
\end{align*}
are selection matrices to select the needed elements, that is, $x_i=R_i\boldsymbol{x}^i$ and $\boldsymbol{x}_{-i}^i=S_i\boldsymbol{x}^i$. $J_i(x_i,\boldsymbol{x}_{-i}^i)$ is the agent $i$'s cost determined by its own strategy and the estimation of strategies of other agents. The extra correction term $R_i\sum_{j\in\mathcal{I}}a_{ij}(\boldsymbol{x}^i-\boldsymbol{x}^j)$ is instrumental in the agreement of agents' estimation vectors.  The strategy $x_i$ is expected to evolve in the direction of any subgradient of cost function $J_i$ to the NE. Thus, the dynamics of auxiliary state $v_i$ can be represented as a differential inclusion. Define $\boldsymbol{x}$ $=$ $\col(\boldsymbol{x}^1,\ldots,\boldsymbol{x}^N)$. A set-valued map $\mathbf{F}(\boldsymbol{x}): \Omega^N\rightrightarrows \mathbb{R}^{Nn}$  is defined by $\mathbf{F}(\boldsymbol{x})$ $=$ $\prod_{i=1}^N\partial_{x_i} J_i(x_i,\boldsymbol{x}_{-i}^i)$. The following assumption is given for the map $\mathbf{F}(\boldsymbol{x})$.

\begin{assumption} \label{as2}
The set-valued map $\mathbf{F}(\boldsymbol{x})$ is Lipschitz continuous with Lipschitz constant $\theta>0$. It also satisfies that $(x-x')^T(\boldsymbol{d}-\boldsymbol{d}')\geq w\|\boldsymbol{x}-\boldsymbol{x}'\|^2$, $w>0$,$\forall \boldsymbol{x}, \boldsymbol{x}'\in \Omega^N$, $\boldsymbol{d}\in \mathbf{F}(\boldsymbol{x}), \boldsymbol{d}'\in \mathbf{F}(\boldsymbol{x}')$.
\end{assumption}

Assumption \ref{as4} ensures that the noncooperative game $G$ has a unique NE \cite[Theorem 2]{Rosen.1965}. The assumptions of Lipschitz continuity and monotonicity of involved maps are also used in \cite{Deng.2019b, Zeng.2019, Gadjov.2019}. Assumption \ref{as2} is an extension of Assumption \ref{as4} from the strategy space $\Omega$ to its augmented space $\Omega^N$. The assumption on the strong monotonicity of $F(x)$ can hold under that each $J_i(x_i,x_{-i})$ is a convex function. For the double-integrator dynamics, the assumption on set-valued map $\mathbf{F}(\boldsymbol{x})$ is stronger than that in \cite{Gadjov.2019, Zeng.2019} and it can guarantee the convergence of the proposed strategy-updating rule.
The following analysis mainly focuses on undirected graphs, which is assumed to be connected in the following assumption.

\begin{assumption} \label{as3}
The undirected communication graph is connected.
\end{assumption}

Let $x=\col(x_1,\ldots,x_N)$, $v=\col(v_1,\ldots,v_N)$, $\boldsymbol{x}$ $=$ $\col(\boldsymbol{x}^1,$ \\$\ldots,\boldsymbol{x}^N)$, $\mathcal{R}$ $=$ $\blk\{R_1,\ldots,R_N\}$, $\mathcal{S}$ $=$ $\blk\{S_1,\ldots,S_N\}$,  $\mathcal{S}\boldsymbol{x}$ $=$ $\col(\boldsymbol{x}_{-1}^1,\ldots,\boldsymbol{x}_{-N}^N)$, and $\boldsymbol{L}$ $=$ $L\otimes I_{Nn}$. A compact form of the designed rule \eqref{csur} for all agents is given as follows,
\begin{equation} \label{cf1}
\begin{aligned}
\dot{x}&=v,\\
\dot{v}&\in-kv-\mathbf{F}(\boldsymbol{x})-\frac{\alpha}{k}\mathcal{R}\boldsymbol{L}\boldsymbol{x},\\
\mathcal{S}\dot{\boldsymbol{x}}&=-\alpha \mathcal{S}\boldsymbol{L}\boldsymbol{x}.
\end{aligned}
\end{equation}

\begin{lemma}
Suppose that Assumption \ref{as1} holds. $x^*$ is the NE of game $G=(\mathcal{I},\Omega,J)$ if and only if $(x^*,\boldsymbol{0}_{Nn}, \boldsymbol{1}_N\otimes x^*)$ is the equilibrium of dynamic system \eqref{cf1}.
\end{lemma}

The proof is straightforward. We omit it for limited space.

\begin{theorem} \label{theorem1}
Suppose that Assumptions \ref{as1}-\ref{as3} hold and the parameters $k$ and $\alpha$ satisfy that $k>\max\{\frac{2\theta}{w},\theta+\sqrt{\theta^2+\alpha\|\mathcal{R}\boldsymbol{L}\|},\frac{\|\mathcal{R}\boldsymbol{L}\|}{\lambda_2}\}$ and  $\alpha(\lambda_2-\frac{\|\mathcal{R}\boldsymbol{L}\|}{k})>\theta+\frac{k^2\theta}{4}$, respectively. Agents with dynamics \eqref{ad} follow the strategy-updating rule \eqref{csur}. Then, all agents' strategies can asymptotically converge to the unique NE of game $G=(\mathcal{I},\Omega,J)$.
\end{theorem}

Proof: Define $\tilde{x}=x-x^*$, $\tilde{v}=v-v^*$, and $\tilde{\boldsymbol{x}}=\boldsymbol{x}-\boldsymbol{1}_{N}\otimes x^*$. The system \eqref{cf1} can be transformed into
\begin{align} \label{bs1}
[\dot{\tilde{x}},\dot{\tilde{v}},\mathcal{S}\dot{\tilde{\boldsymbol{x}}}]^T\in \mathcal{F}(\tilde{x},\tilde{v},\tilde{\boldsymbol{x}}),
\end{align}
where $\mathcal{F}(\tilde{x},\tilde{v},\tilde{\boldsymbol{x}})=\begin{bmatrix} \tilde{v}\\ -k\tilde{v}-\mathbf{F}(\boldsymbol{x})+\mathbf{F}(\boldsymbol{x}^*)-\frac{\alpha}{k}\mathcal{R}\boldsymbol{L}\tilde{\boldsymbol{x}} \\ -\alpha \mathcal{S}\boldsymbol{L}\tilde{\boldsymbol{x}} \end{bmatrix}$.
Recall the definition of $\mathbf{F}(\boldsymbol{x})$. It follows from  Lemma \ref{lemma1} that $\boldsymbol{0}_{Nn}\in \mathbf{F}(\boldsymbol{x}^*)$.

Consider a Lyapunov candidate function $V=\frac{1}{2}(\|\tilde{v}\|^2+\|k\tilde{x}+\tilde{v}\|^2+\tilde{\boldsymbol{x}}^T\mathcal{S}^T\mathcal{S}\tilde{\boldsymbol{x}})$.
The set-valued Lie derivative of $V$ with respect to $\mathcal{F}$ is given by
\begin{align*}
\begin{split}
\mathcal{L}_{\mathcal{F}}V&=\{\zeta\in\mathbb{R}: \zeta=-k\|\tilde{v}\|^2-2\tilde{v}^T(\boldsymbol{d}-\boldsymbol{d}^*)-\frac{2\alpha}{k}\tilde{v}^T\mathcal{R}\boldsymbol{L}\tilde{\boldsymbol{x}}\\
&\ \ \ \  -k\tilde{x}^T(\boldsymbol{d}-\boldsymbol{d}^*)\!-\!\alpha\tilde{\boldsymbol{x}}^T\boldsymbol{L}\tilde{\boldsymbol{x}},\boldsymbol{d}\in \mathbf{F}(\boldsymbol{x}),\boldsymbol{d}^*\in \mathbf{F}(\boldsymbol{x}^*)\}.
\end{split}
\end{align*}

Let $\tilde{\boldsymbol{x}}^{c}=\frac{1}{N}\boldsymbol{1}_N\boldsymbol{1}_N^T\otimes I_{Nn} \tilde{\boldsymbol{x}}$, and $\tilde{\boldsymbol{x}}^{o}=(I_{N^2n}-\frac{1}{N}\boldsymbol{1}_N\boldsymbol{1}_N^T\otimes I_{Nn} ) \tilde{\boldsymbol{x}}$. Then, $\tilde{\boldsymbol{x}}\in \mathbb{R}^{N^2n}$ can be decomposed into two components. One is in the consensus subspace and the other is in the orthogonal complement of the consensus subspace, that is, $\tilde{\boldsymbol{x}}=\tilde{\boldsymbol{x}}^{c}+\tilde{\boldsymbol{x}}^{o}$. Since $\tilde{\boldsymbol{x}}^{c}=\boldsymbol{1}_N\otimes x$ for some $x\in\mathbb{R}^{Nn}$, it follows that $\boldsymbol{L}\tilde{\boldsymbol{x}}^{c}=\boldsymbol{0}_{N^2n}$, and $(\tilde{\boldsymbol{x}}^{o})^T\boldsymbol{L}\tilde{\boldsymbol{x}}^{o}\geq \lambda_2 \|\tilde{\boldsymbol{x}}^{o}\|^2$, where $\lambda_2$ is the second least eigenvalue of $L$.  From the definition of $\tilde{\boldsymbol{x}}^{c}$ and $\tilde{\boldsymbol{x}}^{o}$, $(\tilde{\boldsymbol{x}}^{c})^T\tilde{\boldsymbol{x}}^{o}=0$. Thus, $\|\tilde{\boldsymbol{x}}\|^2=\|\tilde{\boldsymbol{x}}^{c}\|^2+\|\tilde{\boldsymbol{x}}^{o}\|^2$.  Define $\boldsymbol{d}'=\boldsymbol{1}_N\otimes x$ for some $x\in\Omega$. If $\boldsymbol{x}^i=x$, $i\in\mathcal{I}$, $\partial_{x_i} J_i(x_i,\boldsymbol{x}_{-i}^i)=\partial_{x_i}J_i(x_i,x_{-i})$. It follows from the definitions of $F(x)$ and $\mathbf{F}(\boldsymbol{x})$ that $\mathbf{F}(\boldsymbol{1}_N\otimes x)=F(x)$. Under Assumptions \ref{as4} and \ref{as2}, it follows that
\begin{equation}
\begin{aligned}
-2\tilde{v}^T(\boldsymbol{d}-\boldsymbol{d}^*)&=-2\tilde{v}^T(\boldsymbol{d}-\boldsymbol{d}')-2\tilde{v}^T(\boldsymbol{d}'-\boldsymbol{d}^*)\\
& \leq 2\theta\|\tilde{v}\|\|\tilde{\boldsymbol{x}}^{o}\|+2\theta\|\tilde{v}\|\|\tilde{x}\|\\
& \leq 2\theta\|\tilde{v}\|^2+\theta\|\tilde{\boldsymbol{x}}^{o}\|^2+\theta\|\tilde{x}\|,
\end{aligned}
\end{equation}

\begin{equation}
\begin{aligned}
-\frac{2\alpha}{k}\tilde{v}^T\mathcal{R}\boldsymbol{L}\tilde{\boldsymbol{x}}&\leq \frac{2\alpha}{k}\|\mathcal{R}\boldsymbol{L}\|\|\tilde{v}\|\|\tilde{\boldsymbol{x}}^o\|\\
&\leq \frac{\alpha\|\mathcal{R}\boldsymbol{L}\|}{k}\|\tilde{v}\|^2+\frac{\alpha\|\mathcal{R}\boldsymbol{L}\|}{k}\|\|\tilde{\boldsymbol{x}}^o\|^2,
\end{aligned}
\end{equation}
and
\begin{equation}
\begin{aligned}
-k\tilde{x}^T(\boldsymbol{d}\!-\!\boldsymbol{d}^*)&\!\leq\! -k(x\!-\!x^*)^T(\boldsymbol{d}\!-\!\boldsymbol{d}')\!-\!k(x\!-\!x^*)^T(\boldsymbol{d}'\!-\!\boldsymbol{d}^*)\\
&\leq k\theta\|\tilde{x}\|\|\tilde{\boldsymbol{x}}^{o}\|-kw\|\tilde{x}\|^2\\
& \leq \theta\|\tilde{x}\|^2+\frac{k^2\theta}{4}\|\tilde{\boldsymbol{x}}^{o}\|^2-kw\|\tilde{x}\|^2.
\end{aligned}
\end{equation}

Thus,
\begin{align*}
\zeta&\leq-(kw-2\theta)\|\tilde{x}\|^2-(k-2\theta-\frac{\alpha\|\mathcal{R}\boldsymbol{L}\|}{k})\|\tilde{v}\|^2\\
&\ \ \ -(\alpha\lambda_2-\theta-\frac{\alpha\|\mathcal{R}\boldsymbol{L}\|}{k}-\frac{k^2\theta}{4})\|\tilde{\boldsymbol{x}}^o\|^2.
\end{align*}
Since $\zeta$ is arbitrary, it follows that
\begin{align*}
\max \mathcal{L}_{\mathcal{F}} V &\leq -(kw-2\theta)\|\tilde{x}\|^2-(k-2\theta-\frac{\alpha\|\mathcal{R}\boldsymbol{L}\|}{k})\|\tilde{v}\|^2\\
&\ \ \ -(\alpha\lambda_2-\theta-\frac{\alpha\|\mathcal{R}\boldsymbol{L}\|}{k}-\frac{k^2\theta}{4})\|\tilde{\boldsymbol{x}}^o\|^2.
\end{align*}

If $k=\max\{\frac{2\theta}{w},\theta+\sqrt{\theta^2+\alpha\|\mathcal{R}\boldsymbol{L}\|},\frac{\|\mathcal{R}\boldsymbol{L}\|}{\lambda_2}\}$, and $\alpha(\lambda_2-\frac{\|\mathcal{R}\boldsymbol{L}\|}{k})>\theta+\frac{k^2\theta}{4}$, $\max\mathcal{L}_{\mathcal{F}}V<0$ with $\tilde{x}\neq \boldsymbol{0}_{Nn}$, $\tilde{v}\neq \boldsymbol{0}_{Nn}$, or $\tilde{\boldsymbol{x}}\neq \boldsymbol{0}_{N^2n}$. $\max\mathcal{L}_{\mathcal{F}}V=0$ only if $\tilde{x}= \boldsymbol{0}_{Nn}$, $\tilde{v}= \boldsymbol{0}_{Nn}$, and $\tilde{\boldsymbol{x}}=\boldsymbol{0}_{N^2n}$, which indicates that  all agents' strategies arrive at the NE. Recall that $V$ is a continuously differentiable, radially unbounded and positive definite function. It follows from system (6) that the origin is the equilibrium point. According to Corollary 4.2 and Theorem 4.4 in \cite{Khalil.2002}, the largest invariant set is given by
\begin{equation} \label{is}
\begin{split}
M&=\{\tilde{x}\in\mathbb{R}^{Nn},\tilde{v}\in\mathbb{R}^{Nn}, \tilde{\boldsymbol{x}}\in\mathbb{R}^{N^2n}|\tilde{x}=\boldsymbol{0}_{Nn},\\
&\ \ \ \ \ \ \tilde{v}=\boldsymbol{0}_{Nn},\tilde{\boldsymbol{x}}=\boldsymbol{0}_{N^2n}\}.
\end{split}
\end{equation}

By Lemma \ref{lemma2}, any trajectory of \eqref{bs1} starting  from an initial condition $(\tilde{x}_0,\tilde{v}_0,\tilde{\boldsymbol{x}}_0)$ converges to the invariant set $M$.
Thus, $(x,v,\boldsymbol{x})$ converges to the equilibrium $(x^*,\boldsymbol{0}_{Nn},\boldsymbol{1}_N\otimes x^*)$ as $t\rightarrow \infty$.
It indicates that all agents' strategies can reach to the NE of noncooperative game $G=(\mathcal{I},\Omega,J)$. $\hfill\blacksquare$

\begin{remark}
To ensure that each agent estimates the strategies of others accurately, it is necessary to assume the connectivity of communication graphs, which is a global property of communication graphs. In the case that the communication topology is the prior knowledge to agents, parameters $\alpha$ and $k$ can be selected to satisfy the conditions in Theorems \ref{theorem1}-\ref{theorem4}. In addition, $\|\mathcal{R}\boldsymbol{L}\|$ only depends on the maximal degree of nodes in the graph, due to the special structure of $\mathcal{R}$. If the overall structure of communication graph is unknown, the total number of agents is necessary to be known to estimate the algebraic connectivity $\lambda_2$ and $\|\mathcal{R}\boldsymbol{L}\|$.
The conditions involving the Laplacian matrix of the graph can be relaxed by eigenvalue estimations or adaptive gains, such as the adaptive algorithm proposed in \cite{Persis.2019}.
\end{remark}

The above result can be extended to the weight-balanced and strongly connected directed graphs (digraphs). To avoid any confusion, we denote the second smallest eigenvalue of $\frac{1}{2}(L+L^T)$ by $\hat{\lambda}_2$.

\textbf{Corollary 1.} Let $\mathcal{G}$ be a weight-balanced and strongly connected digraph. Under Assumptions \ref{as1} and \ref{as2}, If $k>\max\{\frac{2\theta}{w},\theta+\sqrt{\theta^2+\alpha\|\mathcal{R}\boldsymbol{L}\|},\frac{\|\mathcal{R}\boldsymbol{L}\|}{\hat{\lambda}_2}\}$ and  $\alpha(\hat{\lambda}_2-\frac{\|\mathcal{R}\boldsymbol{L}\|}{k})>\theta+\frac{k^2\theta}{4}$ , all agents' strategies, who have dynamics \eqref{ad} and follow the strategy-updating rule \eqref{csur}, can asymptotically converge to the unique NE of game $G=(\mathcal{I},\Omega,J)$.

Proof:
The proof is similar to that for Theorem \ref{theorem1}. The difference is the treatment on $\frac{1}{2}\tilde{\boldsymbol{x}}^T\mathcal{S}^T\mathcal{S}\tilde{\boldsymbol{x}}$ in the Lyapunov function $V$. For a weighted-balanced and strongly connected digraph, the set-valued Lie derivative of $V$ with respect to $\mathcal{F}$ is
\begin{align*}
\begin{split}
\mathcal{L}_{\mathcal{F}}V&=\{\zeta\in\mathbb{R}:\zeta=-k\|\tilde{v}\|^2-2\tilde{v}^T(\boldsymbol{d}-\boldsymbol{d}^*)-\frac{2\alpha}{k}\tilde{v}^T\mathcal{R}\boldsymbol{L}\tilde{\boldsymbol{x}}\\
&\ \ \ \  -k\tilde{x}^T(\boldsymbol{d}-\boldsymbol{d}^*)-\frac{\alpha}{2}\tilde{\boldsymbol{x}}^T(\boldsymbol{L}+\boldsymbol{L}^T)\tilde{\boldsymbol{x}},\\
&\ \ \ \ \boldsymbol{d}\in \mathbf{F}(\boldsymbol{x}),\boldsymbol{d}^*\in \mathbf{F}(\boldsymbol{x}^*)\}.
\end{split}
\end{align*}

The rest analysis is similar to that in the proof of Theorem \ref{theorem1} and omitted for saving the space. $\hfill\blacksquare$

\section{Distributed strategy-updating rule with discrete-time communication}
In this section, discrete-time communication schemes for strategy-updating rule \eqref{csur} are explored. The implementation of strategy-updating rule \eqref{csur} requires agents to communicate each other in continuous time, which facilitates the theoretical analysis. Considering the cost and execution mechanism of communication in practical scenarios, we study the strategy-updating rule in discrete-time communication schemes. In this section, communication topologies described by undirected graphs are considered.

Let $\{t_k^i\}_{k=1}^\infty \in \mathbb{R}_{\geq 0}$, such that $t_k^i<t_{k+1}^i$, denote the time sequence at which agent $i$ broadcasts its  estimation state $\boldsymbol{x}^i(t_k^i)$ to its neighbors, for all $i\in\mathcal{I}$. Before the next time $t_{k+1}^i$, $\hat{\boldsymbol{x}}^i(t)=\boldsymbol{x}^i(t_k^i)$ for $t\in[t_k^i,t_{k+1}^i)$.
Sometimes $t$ is omitted for simplicity. For agent $i$, the strategy-updating rule \eqref{csur} with discrete-time communication is given by
\begin{equation} \label{dsur}
\begin{aligned}
\dot{x}_i&=v_i,\\
\dot{v}_i&\in-kv_i-\partial_{x_i}J_i(x_i,\boldsymbol{x}_{-i}^i)-\frac{\alpha}{k}R_i\sum_{j\in\mathcal{I}}a_{ij}(\hat{\boldsymbol{x}}^i-\hat{\boldsymbol{x}}^j),\\
\dot{\boldsymbol{x}}_{-i}^i&=-\alpha S_i\sum_{j\in \mathcal{I}}a_{ij}(\hat{\boldsymbol{x}}^i-\hat{\boldsymbol{x}}^j).
\end{aligned}
\end{equation}

Next, two discrete-time communication schemes are proposed for agents to interact with each other at discrete-time instants. Under these schemes, it is analyzed that all agents' strategies can converge asymptotically to the NE of game $G=(\mathcal{I},\Omega,J)$. One is a periodic communication scheme and the other is an event-triggered communication scheme.

\subsection{Periodic Communications}
In the periodic communication scheme, all agents communicate with each other synchronously at time interval $\triangle$, i.e., sampling period $\triangle =t_{k+1}^i-t_k^i$ for all $i\in\mathcal{I}$. Theorem \ref{theorem2} presents an upper bound on the size of execution cycle of communications among agents over an undirected graph.

\begin{theorem} \label{theorem2}
Suppose that Assumptions \ref{as1}-\ref{as3} hold and the parameters $k$ and $\alpha$ satisfies that $k>\max\{\frac{2\theta}{w},\theta+\sqrt{\theta^2+2\alpha\|\mathcal{R}\boldsymbol{L}\|},\frac{\|\mathcal{R}\boldsymbol{L}\|}{\lambda_2}+\frac{\|\boldsymbol{L}\|}{2\lambda_2}\}$ and $\alpha(\lambda_2-\frac{\|\mathcal{R}\boldsymbol{L}\|}{k}-\frac{\|\boldsymbol{L}\|}{2k})>\theta+\frac{k^2\theta}{4}$, respectively. Each agent communicates synchronously with its neighbors over the graph $\mathcal{G}$ every $\triangle$ seconds, starting at zero, where $\triangle\in(0,\tau)$, $\tau$ is the upper bound of communication intervals given by
\begin{align} \label{tau}
\tau=\frac{1}{a}\ln(1+\frac{a\xi}{a+b+b\xi}).
\end{align}
In \eqref{tau},  $\xi^2=\frac{(\alpha \lambda_2-\theta)k-\alpha\|\mathcal{R}\boldsymbol{L}\|-k^3\theta/4-\alpha\|\boldsymbol{L}\|/2}{\alpha\|\mathcal{R}\boldsymbol{L}\|+\alpha \|\boldsymbol{L}\|k^2/2}$, $a=\frac{\theta+1}{\alpha}$, and   $b=\alpha\|\mathcal{S}^T\mathcal{S}\boldsymbol{L}\|$. Agents with dynamics \eqref{ad} follow the strategy-updating rule \eqref{dsur}.
Then, all agents' strategies can asymptotically reach to the unique NE of game $G=(\mathcal{I},\Omega,J)$.
\end{theorem}
Proof:
First, similar to the analysis of Theorem 1, we transfer the equilibrium to the origin, which is similar in Theorem \ref{theorem1}. Let $\boldsymbol{e}_i=\hat{\boldsymbol{x}}^i(t_k^i)-\boldsymbol{x}^i(t), \forall i\in \mathcal{I}$ and $\boldsymbol{e}=\col(\boldsymbol{e}_1,\ldots,\boldsymbol{e}_N)$ for $\boldsymbol{e}=\hat{\boldsymbol{x}}-\boldsymbol{x}$. \eqref{dsur} can be written in a compact form
\begin{align}\label{did}
[\dot{\tilde{x}},\dot{\tilde{v}}, \mathcal{S}\dot{\tilde{\boldsymbol{x}}}]^T\in \bar{\mathcal{F}}(\tilde{x},\tilde{v},\tilde{\boldsymbol{x}}),
\end{align}
where  $\bar{\mathcal{F}}(\tilde{x},\tilde{v},\tilde{\boldsymbol{x}})=\begin{bmatrix} \tilde{v}\\ -k\tilde{v}-\mathbf{F}(\boldsymbol{x})+\mathbf{F}(\boldsymbol{x}^*)-\frac{\alpha}{k}\mathcal{R}\boldsymbol{L}(\boldsymbol{e}+\tilde{\boldsymbol{x}}) \\ -\alpha \mathcal{S}\boldsymbol{L}(\boldsymbol{e}+\tilde{\boldsymbol{x}}) \end{bmatrix}$.

Recall the definition of Lyapunov function $V$ in Theorem \ref{theorem1}.  The set-valued Lie derivative of $V$ with respect to $\bar{\mathcal{F}}$ is given by
\begin{equation} \label{lv1}
\begin{split}
\mathcal{L}_{\bar{\mathcal{F}}}V&=\{\zeta\in\mathbb{R}:\zeta=-k\|\tilde{v}\|^2-2\tilde{v}^T(\boldsymbol{d}-\boldsymbol{d}^*)-\frac{2\alpha}{k}\tilde{v}^T\mathcal{R}\boldsymbol{L}\tilde{\boldsymbol{x}}\\
&\ \ \ \  -k\tilde{x}^T(\boldsymbol{d}-\boldsymbol{d}^*)\!-\!\alpha\tilde{\boldsymbol{x}}^T\boldsymbol{L}\tilde{\boldsymbol{x}} -\frac{2\alpha}{k}\tilde{v}^T\mathcal{R}\boldsymbol{L}\boldsymbol{e}-\alpha\tilde{\boldsymbol{x}}^T\boldsymbol{L}\boldsymbol{e},\\
&\ \ \ \ \boldsymbol{d}\in \mathbf{F}(\boldsymbol{x}),\boldsymbol{d}^*\in \mathbf{F}(\boldsymbol{x}^*)\}.
\end{split}
\end{equation}
Similar to the analysis in Theorem \ref{theorem1}, it follows that
\begin{equation} \label{zeta2}
\begin{split}
\zeta&\leq -(kw-2\theta)\|\tilde{x}\|^2-(k-2\theta-\frac{\alpha\|\mathcal{R}\boldsymbol{L}\|}{k})\|\tilde{v}\|^2\\
&\ \ \ -(\alpha\lambda_2-\theta-\frac{\alpha\|\mathcal{R}\boldsymbol{L}\|}{k}-\frac{k^2\theta}{4})\|\tilde{\boldsymbol{x}}^o\|^2\\
&\ \ \ +\frac{2\alpha}{k}\|\mathcal{R}\boldsymbol{L}\|\|\tilde{v}\|\|\boldsymbol{e}\|+\alpha\|\boldsymbol{L}\|\|\tilde{\boldsymbol{x}}^o\|\|\boldsymbol{e}\|
\end{split}
\end{equation}

By Young Inequality, we have that
\begin{equation} \label{16}
\frac{2\alpha}{k}\|\mathcal{R}\boldsymbol{L}\|\|\tilde{v}\|\|\boldsymbol{e}\|\leq \frac{\alpha\|\mathcal{R}\boldsymbol{L}\|}{k}\|\tilde{v}\|^2+\frac{\alpha\|\mathcal{R}\boldsymbol{L}\|}{k}\|\boldsymbol{e}\|^2,
\end{equation}
and
\begin{equation} \label{17}
\alpha\|\boldsymbol{L}\|\|\tilde{\boldsymbol{x}}^o\|\|\boldsymbol{e}\|\leq \frac{\alpha\|\boldsymbol{L}\|}{2k}\|\tilde{\boldsymbol{x}}^o\|^2+\frac{\alpha k\|\boldsymbol{L}\|}{2}\|\boldsymbol{e}\|.
\end{equation}

Substituting \eqref{16} and \eqref{17} into \eqref{zeta2} yields that
\begin{equation}
\begin{split}
\zeta&\leq-(kw-2\theta)\|\tilde{x}\|^2-(k-2\theta-\frac{2\alpha\|\mathcal{R}\boldsymbol{L}\|}{k})\|\tilde{v}\|^2\\
&\ \ \ -(\alpha\lambda_2-\theta-\frac{\alpha\|\mathcal{R}\boldsymbol{L}\|}{k}-\frac{\alpha\|\boldsymbol{L}\|}{2k}-\frac{k^2\theta}{4})\|\tilde{\boldsymbol{x}}^o\|^2\\
&\ \ \ +(\frac{\alpha\|\mathcal{R}\boldsymbol{L}\|}{k}+\frac{\alpha k\|\boldsymbol{L}\|}{2})\|\boldsymbol{e}\|.
\end{split}
\end{equation}

Since $\zeta$ is arbitrary, we have that
\begin{align*}
\begin{split}
\max\mathcal{L}_{\bar{\mathcal{F}}}V&\leq-(kw-2\theta)\|\tilde{x}\|^2-(k-2\theta-\frac{2\alpha\|\mathcal{R}\boldsymbol{L}\|}{k})\|\tilde{v}\|^2\\
&\ \ \ -(\alpha\lambda_2-\theta-\frac{\alpha\|\mathcal{R}\boldsymbol{L}\|}{k}-\frac{\alpha\|\boldsymbol{L}\|}{2k}-\frac{k^2\theta}{4})\|\tilde{\boldsymbol{x}}^o\|^2\\
&\ \ \ +(\frac{\alpha\|\mathcal{R}\boldsymbol{L}\|}{k}+\frac{\alpha k\|\boldsymbol{L}\|}{2})\|\boldsymbol{e}\|,
\end{split}
\end{align*}
where $k=\max\{\frac{2\theta}{w},\theta+\sqrt{\theta^2+2\alpha\|\mathcal{R}\boldsymbol{L}\|},\frac{\|\mathcal{R}\boldsymbol{L}\|}{\lambda_2}+\frac{\|\boldsymbol{L}\|}{2\lambda_2}\}$ and $\alpha(\lambda_2-\frac{\|\mathcal{R}\boldsymbol{L}\|}{k}-\frac{\|\boldsymbol{L}\|}{2k})>\theta+\frac{k^2\theta}{4}$.
Let $\xi^2=\frac{(\alpha \lambda_2-\theta)k-\alpha\|\mathcal{R}\boldsymbol{L}\|-k^3\theta/4-\alpha\|\boldsymbol{L}\|/2}{\alpha\|\mathcal{R}\boldsymbol{L}\|+\alpha \|\boldsymbol{L}\|k^2/2}$. If $\|\boldsymbol{e}\|^2<\xi^2 \|\tilde{\boldsymbol{x}}^{o}(t)\|^2<\xi^2 \|\tilde{\boldsymbol{x}}(t)\|^2$, $t\in[t_k,t_{k+1})$, $\max\mathcal{L}_{\bar{\mathcal{F}}}V_1<0$ for all $t>0$. It is clear that at each communication time $t_k$, $\|\boldsymbol{e}\|=0$. Then, $\boldsymbol{e}(t)$ grows until next communication time $t_{k+1}$ and becomes zero again. The following analysis shows the upper bound of the communication intervals by examining the time period it takes for $q=\|\boldsymbol{e}\|/\|\tilde{\boldsymbol{x}}(t)\|$ to evolve from zero to $\xi$.

\begin{equation}
\begin{aligned}
\dot{q}&=\frac{\boldsymbol{e}^T\dot{\boldsymbol{e}}}{\|\boldsymbol{e}\|\|\tilde{\boldsymbol{x}}\|}-\frac{\|\boldsymbol{e}\|\tilde{\boldsymbol{x}}^T\dot{\tilde{\boldsymbol{x}}}}{\|\tilde{\boldsymbol{x}}\|^3}\leq \frac{\|\dot{\boldsymbol{e}}\|}{\|\tilde{\boldsymbol{x}}\|}+\frac{\|\boldsymbol{e}\|\|\dot{\tilde{\boldsymbol{x}}}\|}{\|\tilde{\boldsymbol{x}}\|^2}\\
&\leq (1+q)\frac{\|\dot{\tilde{\boldsymbol{x}}}\|}{\|\tilde{\boldsymbol{x}}\|}.
\end{aligned}
\end{equation}
The second inequality follows from the definition of $q$ and the fact that $\|\dot{\boldsymbol{e}}\|\leq \|\dot{\tilde{\boldsymbol{x}}}\|$.

In addition,
\begin{align*}
\frac{\|\dot{\tilde{\boldsymbol{x}}}\|}{\|\tilde{\boldsymbol{x}}\|}&=\frac{\tilde{\boldsymbol{x}}^T\dot{\tilde{\boldsymbol{x}}}}{\|\tilde{\boldsymbol{x}}\|^2}
= \frac{\tilde{x}^T\tilde{v}-\alpha\tilde{\boldsymbol{x}}^T\mathcal{S}^T\mathcal{S}\boldsymbol{L}\boldsymbol{e}-\alpha\tilde{\boldsymbol{x}}^T\mathcal{S}^T\mathcal{S}\boldsymbol{L}\tilde{\boldsymbol{x}}}{\|\tilde{\boldsymbol{x}}\|^2}\\
&\leq \frac{\|\mathcal{R}\tilde{\boldsymbol{x}}\|\|\tilde{v}\|}{\|\tilde{\boldsymbol{x}}\|^2}+\alpha\frac{\|\mathcal{S}^T\mathcal{S}\boldsymbol{L}\|\|\boldsymbol{e}\|}{\|\tilde{\boldsymbol{x}}\|}+\alpha\|\mathcal{S}^T\mathcal{S}\boldsymbol{L}\|\\
&\leq \frac{\|\tilde{v}\|}{\|\tilde{\boldsymbol{x}}\|}+\alpha\|\mathcal{S}^T\mathcal{S}\boldsymbol{L}\|(1+q).
\end{align*}

It follows from the evolution of \eqref{dsur} that $\frac{\|\tilde{v}\|}{\|\tilde{\boldsymbol{x}}\|}$ $\leq$ $\int_0^t e^{(\tau-t)}\frac{\|F(\boldsymbol{x})-F(\boldsymbol{x}^*)+\mathcal{R}\boldsymbol{L}(\boldsymbol{e}+\tilde{\boldsymbol{x}}(\tau))\|}{\alpha\|\mathcal{S}\boldsymbol{L}(\boldsymbol{e}+\tilde{\boldsymbol{x}}(\tau))\|}d\tau$ with $v(0)=0$. By the integration mean value theorem, $\frac{\|\tilde{v}\|}{\|\tilde{\boldsymbol{x}}\|}\leq\frac{(1-e^{-t})}{\alpha\|\mathcal{S}\boldsymbol{L}(\boldsymbol{e}+\tilde{\boldsymbol{x}}(s))\|}\times$ $ (\theta\|\tilde{\boldsymbol{x}}(s)\|+\|\mathcal{R}\boldsymbol{L}(\boldsymbol{e}+\tilde{\boldsymbol{x}}(s))\|)$ for some fixed $s\in (0,t)$. It yields that $\frac{\|\tilde{v}\|}{\|\tilde{\boldsymbol{x}}\|}$ $\leq$ $\frac{\theta+1}{\alpha}$ by the fact that $\|\tilde{\boldsymbol{x}}(s)\|$ $<$ $\|\mathcal{S}\boldsymbol{L}(\boldsymbol{e}+\tilde{\boldsymbol{x}}(s))\|$ and $\|\mathcal{R}\boldsymbol{L}(\boldsymbol{e}+\tilde{\boldsymbol{x}}(s))\|$ $<$ $\|\mathcal{S}\boldsymbol{L}(\boldsymbol{e}+\tilde{\boldsymbol{x}}(s))\|$.
 Thus,
\begin{equation}
\dot{q}\leq (1+q)(\frac{\theta+1}{\alpha}+\alpha\|\mathcal{S}^T\mathcal{S}\boldsymbol{L}\|(1+q)).
\end{equation}

Using the Comparison Lemma in \cite{Khalil.2002}, we have that $q(t,q_0)\leq \varphi(t,\varphi_0)$, where $\varphi(t,\varphi_0)$ is the solution of $\dot{\varphi}=\frac{\theta+1}{\alpha}(1+\varphi)+\alpha\|\mathcal{S}^T\mathcal{S}\boldsymbol{L}\|(1+\varphi)^2$ with initial state $\varphi(0,\varphi_0)=\varphi_0$. Then,
\begin{align*}
q(t,0)\leq \varphi(t,0)=\frac{(a+b)(e^{at}-1)}{a+b(1-e^{at})},
\end{align*}
where $a=\frac{\theta+1}{\alpha}$ and $b=\alpha\|\mathcal{S}^T\mathcal{S}\boldsymbol{L}\|$.

The time $\tau$ when $\varphi(\tau,0)=\xi$ is given by
\begin{equation*}
\tau=\frac{1}{a}\ln(1+\frac{a\xi}{a+b+b\xi}).
\end{equation*}
Then, for $\{t_{k+1}-t_k\}<\tau$, $\|\boldsymbol{e}\|<\xi \|\tilde{\boldsymbol{x}}(t)\|$.  Thus, $\max\mathcal{L}_{\bar{\mathcal{F}}}V=0 $ if $\tilde{x}= \boldsymbol{0}_{Nn}$, $\tilde{v}= \boldsymbol{0}_{Nn}$ and $\tilde{\boldsymbol{x}}=\boldsymbol{0}_{N^2n}$; and $\max\mathcal{L}_{\bar{\mathcal{F}}}V\neq0$, otherwise. The largest invariant set is the same as \eqref{is}. It follows from Lemma \ref{lemma2} that the system \eqref{did} asymptotically converges to the origin, which indicates that all agents' strategies can reach to the NE of noncooperative game $G$. $\hfill\blacksquare$

\begin{remark}
The communication period determined by $\tau$ in Theorem \ref{theorem2} relies on the communication graph, cost functions of players, and the designed parameter $\alpha$. When cost functions are given, and the graph and parameter $\alpha$ are fixed, $\tau$ can be determined by \eqref{tau}.
\end{remark}

\subsection{Dynamic Event-triggered Communications}
Although the periodic communications can be realized easily, it may degrade the system performance and use communication resources with low efficiency. In the following, an event-triggered communication scheme is designed to overcome these weaknesses.

A dynamic event-triggered mechanism, which was proposed in \cite{Girard.2015}, is utilized here. Intoduce an internal dynamic variable $\eta_i\in\mathbb{R}$ for each agent $i\in \mathcal{I}$, and $\eta_i$ is governed by the following dynamics.
\begin{equation} \label{eta}
\dot{\eta}_i=-b\eta_i+\frac{1}{2} \sum_{j=1}^Na_{ij}\|\hat{\boldsymbol{x}}^i-\hat{\boldsymbol{x}}^j\|^2-(2d_i+\beta_1+\beta_2)\|\hat{\boldsymbol{x}}^i-\boldsymbol{x}^i\|^2,
\end{equation}
where $b\!>\!0$, $\beta_1=\frac{\alpha \|\mathcal{R}\boldsymbol{L}\|}{k}$, $\beta_2=\frac{(\alpha-1)k\|\boldsymbol{L}\|}{2}$, $k=\max\{\frac{2\theta}{w},\theta+\sqrt{\theta^2+2\alpha\|\mathcal{R}\boldsymbol{L}\|},\frac{\|\mathcal{R}\boldsymbol{L}\|}{\lambda_2}+\frac{\|\boldsymbol{L}\|}{2\lambda_2}\}$ and $(\alpha-1)(\lambda_2-\frac{\|\mathcal{R}\boldsymbol{L}\|}{k}-\frac{\|\boldsymbol{L}\|}{2k})>\theta+\frac{k^2\theta}{4}-\frac{\|\mathcal{R}\boldsymbol{L}\|}{k}$.

\begin{theorem} \label{theorem4}
Suppose that Assumptions \ref{as1}-\ref{as3} hold.  Agent $i$ asynchronously communicates with its neighbors over graph $\mathcal{G}$ at times $\{t_k^i\}_{k\in \mathbb{Z}_{>0}}$, starting at $t_0^i=0$, for all $i\in \mathcal{I}$, according to the following dynamic event-triggering rule
\begin{equation} \label{der}
\begin{aligned}
t_{k+1}^i&=\inf\big\{t\in(t_k^i,\infty)|(\beta_1+\beta_2+2d_i)\|\hat{\boldsymbol{x}}^i-\boldsymbol{x}^i\|^2\\
&\ \ \ \ \ \geq\frac{1}{2}\sum_{i=1}^Na_{ij}\|\hat{\boldsymbol{x}}^i-\hat{\boldsymbol{x}}^j\|^2+\rho\eta_i\big\},
\end{aligned}
\end{equation}
where $\beta_1$, $\beta_2$, $k$ and $\alpha$ are defined in \eqref{eta}, and $\rho>0$. Agents with dynamics \eqref{ad} follow the strategy-updating rule \eqref{dsur}.
Then, all agents' strategies can asymptotically evolve to the unique NE of noncooperative game $G=(\mathcal{I},\Omega,J)$.
\end{theorem}

Proof: Consider the Lyapunov function $V$ defined in Theorem \ref{theorem1}, whose set-valued Lie derivative with respect to $\bar{\mathcal{F}}$ is given by \eqref{lv1}. According to the analysis in Theorem \ref{theorem2}, we have that
\begin{align*}
\zeta&=-k\|\tilde{v}\|^2-2\tilde{v}^T(\boldsymbol{d}-\boldsymbol{d}^*)-\frac{2\alpha}{k}\tilde{v}^T\mathcal{R}\boldsymbol{L}\tilde{\boldsymbol{x}}  -k\tilde{x}^T(\boldsymbol{d}\!-\!\boldsymbol{d}^*)\\ &\ \ \ -(\alpha-1)\tilde{\boldsymbol{x}}^T\boldsymbol{L}\tilde{\boldsymbol{x}} -\frac{2\alpha}{k}\tilde{v}^T\mathcal{R}\boldsymbol{L}\boldsymbol{e}-(\alpha-2)\tilde{\boldsymbol{x}}^T\boldsymbol{L}\boldsymbol{e}+s,\\
  &\ \ \ \ \forall \boldsymbol{d}\in \mathbf{F}(\boldsymbol{x}),\boldsymbol{d}^*\in \mathbf{F}(\boldsymbol{x}^*),
\end{align*}
where $s=-\tilde{\boldsymbol{x}}^T\boldsymbol{L}\tilde{\boldsymbol{x}}-2\tilde{\boldsymbol{x}}^T\boldsymbol{L}\boldsymbol{e}=-\hat{\boldsymbol{x}}^T\boldsymbol{L}\hat{\boldsymbol{x}}+\boldsymbol{e}^T\boldsymbol{L}\boldsymbol{e}$.

Similar to the proof of Theorem \ref{theorem2}, the analysis is given as follows. From  $\tilde{\boldsymbol{x}}=\tilde{\boldsymbol{x}}^{c}+\tilde{\boldsymbol{x}}^{o}$ and $\|\tilde{\boldsymbol{x}}\|^2= \|\tilde{\boldsymbol{x}}^{c}\|^2+\|\tilde{\boldsymbol{x}}^{o}\|^2$, it yields that
\begin{align*}
\zeta \leq
&-(kw-2\theta)\|\tilde{x}\|^2-(k-2\theta-\frac{2\alpha\|\mathcal{R}\boldsymbol{L}\|}{k})\|\tilde{v}\|^2\\
&\ \ \ -(\alpha\lambda_2-\theta-\frac{\alpha\|\mathcal{R}\boldsymbol{L}\|}{k}-\frac{k^2\theta}{4}-\frac{(\alpha-1)\|\boldsymbol{L}\|}{2k})\|\tilde{\boldsymbol{x}}^o\|^2\\
&\ \ \ +(\beta_1+\beta_2)\|\boldsymbol{e}\|+s,
\end{align*}
where $\beta_1$ and $\beta_2$ are defined in \eqref{eta}.

From $L=D-A$ and $D+A\geq0$ with the degree matrix $D$ and the adjacent matrix $A$ of graph $\mathcal{G}$, it follows that $\boldsymbol{e}^T\boldsymbol{L}\boldsymbol{e}\leq2\boldsymbol{e}^T(D\otimes I_{Nn})\boldsymbol{e}=2\sum_{i=1}^Nd_i\|\boldsymbol{e}_i\|^2$. Therefore, we have that $s=\frac{1}{2}\sum_{i=1}^N\big(4d_i\|\boldsymbol{e}_i\|^2-\sum_{j=1}^Na_{ij}\|\hat{\boldsymbol{x}}^i-\hat{\boldsymbol{x}}^j\|^2\big)$. Then,
\begin{align*}
\zeta&\leq-(kw-2\theta)\|\tilde{x}\|^2-(k-2\theta-\frac{2\alpha\|\mathcal{R}\boldsymbol{L}\|}{k})\|\tilde{v}\|^2\\
&\ \ \ -(\alpha\lambda_2-\theta-\frac{\alpha\|\mathcal{R}\boldsymbol{L}\|}{k}-\frac{k^2\theta}{4}-\frac{(\alpha-1)\|\boldsymbol{L}\|}{2k})\|\tilde{\boldsymbol{x}}^o\|^2\\
&\ \ \ -\frac{1}{2}\sum_{i=1}^N\sum_{j=1}^Na_{ij}\|\hat{\boldsymbol{x}}^i-\hat{\boldsymbol{x}}^j\|^2+\sum_{i=1}^N(2d_i+\beta_1+\beta_2)\|\boldsymbol{e}_i\|^2.
\end{align*}

Let
\begin{align*}
\bar{\mathcal{F}}_1(\tilde{x},\tilde{v},\tilde{\boldsymbol{x}},\eta_i)\!=\!
\begin{bmatrix}
\begin{smallmatrix}
\tilde{v}\\ -k\tilde{v}-\mathbf{F}(\boldsymbol{x})+\mathbf{F}(\boldsymbol{x}^*)-\frac{\alpha}{k}\mathcal{R}\boldsymbol{L}(\boldsymbol{e}+\tilde{\boldsymbol{x}}) \\ -\alpha \mathcal{S}\boldsymbol{L}(\boldsymbol{e}+\tilde{\boldsymbol{x}}) \\ \!-\!b\eta_i\!+\!\frac{1}{2} \sum_{j=1}^Na_{ij}\|\hat{\boldsymbol{x}}^i\!-\!\hat{\boldsymbol{x}}^j\|^2\!-\!(2d_i\!+\!\beta_1\!+\!\beta_2)\|\hat{\boldsymbol{x}}^i\!-\!\boldsymbol{x}^i\|^2 \end{smallmatrix}
\end{bmatrix}
\end{align*}
be a set-valued map. Consider the Lyapunov candidate function $V_1=V+\sum_{i=1}^N \eta_i(t)$.
Then, \eqref{dsur} and \eqref{eta} are written as
\begin{equation} \label{did2}
[\dot{\tilde{x}},\dot{\tilde{v}}, \mathcal{S}\dot{\tilde{\boldsymbol{x}}}, \dot{\eta_i}]^T\in \bar{\mathcal{F}}_1(\tilde{x},\tilde{v},\tilde{\boldsymbol{x}},\eta_i),
\end{equation}

An upper bound of set-valued Lie derivative of $V_1$ with respect to $\bar{\mathcal{F}}$ is estimated as follows
\begin{equation}
\begin{aligned}
\max\mathcal{L}_{\bar{\mathcal{F}}_1}{V_1}&\leq -(kw-2\theta)\|\tilde{x}\|^2-(k-2\theta-\frac{2\alpha\|\mathcal{R}\boldsymbol{L}\|}{k})\|\tilde{v}\|^2\\
&\ \ \ -(\alpha\lambda_2-\theta-\frac{\alpha\|\mathcal{R}\boldsymbol{L}\|}{k}-\frac{k^2\theta}{4}\\
&\ \ \ -\frac{(\alpha-1)\|\boldsymbol{L}\|}{2k})\|\tilde{\boldsymbol{x}}^o\|^2-b\sum_{i=1}^N\eta_i.
\end{aligned}
\end{equation}

For $t$ $\in$ $[t_k^i,t_{k+1}^i)$, substituting the triggering condition \eqref{der} into the dynamics \eqref{eta} yields that $\dot{\eta}_i$ $\geq$ $-(b+\rho)\eta_i$. Thus, $\eta_i(t)$ $\geq$ $\eta_i(0)e^{-(b+\rho)t}$ for $\eta_i(0)>0$. Therefore, $\max\mathcal{L}_{\bar{\mathcal{F}}_1}{V_1}=0$ if $\tilde{x}= \boldsymbol{0}_{Nn}$, $\tilde{v}= \boldsymbol{0}_{Nn}$, and $\tilde{\boldsymbol{x}}^{o}=\boldsymbol{0}_{N^2n}$; and $\max\mathcal{L}_{\bar{\mathcal{F}}_1}{V_1}\neq0$ otherwise. Similar to the analysis in Theorem \ref{theorem2}, system \eqref{did2} can converge asymptotically to the origin.

Then, we analyze the Zeno behavior by computing a positive lower bound on the event-interval times in the event-triggered process. The lower bound is denoted by $\tau_i \in \mathbb{R}_{\geq 0}$, which is the elapse of the time that $(\beta_1+\beta_2+2d_i)\|\hat{\boldsymbol{x}}^i-\boldsymbol{x}^i(t)\|^2$ evolves from $0$ to $\rho\eta_i$ for all $i\in \mathcal{I}$. Let
\begin{align*}
\varphi=\frac{m\|\hat{\boldsymbol{x}}^i-\boldsymbol{x}^i(t)\|}{\sqrt{\rho\eta_i}},
\end{align*}
where $m=\sqrt{\beta_1+\beta_2+2d_i}$. The derivative of $\varphi$ with respect to $t$ is give by
\begin{align*}
\dot{\varphi}&=-\frac{m(\hat{\boldsymbol{x}}^i-\boldsymbol{x}^i(t))^T\dot{\boldsymbol{x}}^i(t)}{\|\hat{\boldsymbol{x}}^i-\boldsymbol{x}^i(t)\|\sqrt{\rho\eta_i}}-\frac{m\|\hat{\boldsymbol{x}}^i-\boldsymbol{x}^i(t)\|}{2(\rho\eta_i)^{3/2}}\dot{\eta_i}\\
&\leq\frac{m}{\sqrt{\rho\eta_i}}\|\dot{\boldsymbol{x}}^i(t)\|+\frac{b+\rho}{2\rho}\varphi.
\end{align*}
For $t$ $\in$ $[t_k^i,t_{k+1}^i)$, we have that $\frac{m}{\sqrt{\rho\eta_i}}< \infty$ and $\|\dot{\boldsymbol{x}}^i(t)\|$ is bound. Thus, it yields that $\frac{m}{\sqrt{\rho\eta_i}}\|\dot{\boldsymbol{x}}^i(t)\|\leq D$ for some positive constant $D$. It yields that
\begin{align*}
\dot{\varphi}\leq D+\frac{b+\rho}{2\rho}\varphi.
\end{align*}
We have $\varphi(t)\leq \frac{2\rho D}{b+\rho}\big(e^{\frac{b+\rho}{2\rho}(t-t_k^i)}-1\big), \ t\geq t_k^i$ by using the Comparison Lemma \cite[Lemma 3.4]{Khalil.2002} and by the fact that $\|\hat{\boldsymbol{x}}^i-\boldsymbol{x}^i(t_k^i)\|=0$.

Then,
\begin{align*}
\tau_i=\frac{2\rho}{b+\rho}\ln(\frac{b+\rho}{2\rho D}+1).
\end{align*}
which indicates that the Zeno-behavior is excluded in the designed event-triggered scheme \eqref{der}.

$\hfill\blacksquare$

\section{Simulations}
Here, an example on networks of Cournot competition is given to illustrate the effectiveness of the designed continuous-time strategy-updating rule \eqref{csur} and the rule with discrete-time communication \eqref{dsur}, respectively.

The competition among distributed energy resources is considered here, where turbine-generator systems can be described by double integrator agents who communicate with each other on a circle graph \cite{Deng.2019,Deng.2019b}. The cost function of agent $i$ ($i\in \{1,\ldots,5\}$), is
\begin{equation*}
J_i(x_i,x_{-i})=\delta_i+\beta_i|x_i-c_i|+\gamma_ix_i^2-(p-a\sum_{i=1}^Nx_i^2)x_i,
\end{equation*}
where $\delta=[5, 8, 6, 9, 7]^T$, $\beta=[12, 15, 8, 11, 13]^T$, $\gamma=[0.4, 0.5, 0.5, 0.3, 0.3]^T$, $c=[25, 48, 15, 30, 45]^T$, $x(0)=[25, 30, 20, 30, 35]^T$, $p=10$, and $a=0.001$. The cost functions satisfy the Assumptions \ref{as1} and \ref{as2}.

\begin{figure}[!t]
  \centering
  \includegraphics[width=7.5cm]{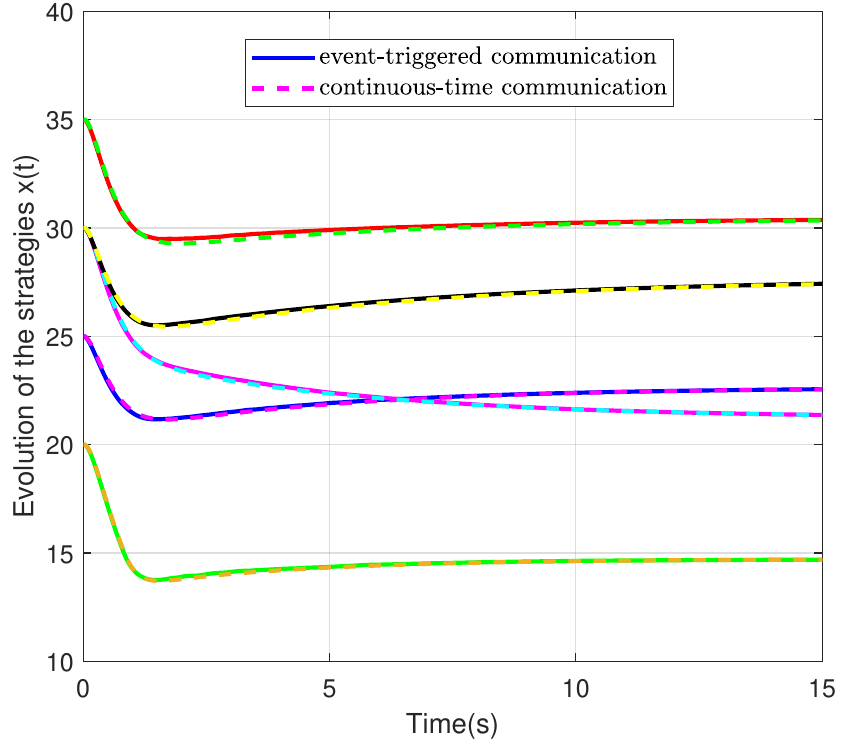}
  \caption{The evolution of the strategies of all agents following \eqref{csur} with continuous-time communications and following \eqref{dsur} with event-triggered communications \eqref{der}}
  \label{fig1}
\end{figure}

Based on the given cost functions and the communication graph, we obtained that $w=0.601$, $\theta=1.001$, and $\lambda_2=1.382$. The parameters $\alpha$, $k$, $b$, and $\rho$ are selected to satisfy the conditions given in Theorems \ref{theorem1}-\ref{theorem4}. In Fig.~\ref{fig1}, the dash lines depict the evolution of all agents' strategies, who follow the continuous-time strategy-updating rule \eqref{csur} with parameters $k=4$ and $\alpha=5$, and the solid lines draw the evolution of all agents' strategies, who follow \eqref{dsur} with the communication scheme based on the event-triggering rule with $k=4$, $\alpha=5$, $b=0.01$, and $\rho=3$. It is seen that all agents' strategies asymptotically converge to the NE of game $G$. Moreover, the strategy-updating rule based on the event-triggered scheme \eqref{der} has a similar convergence performance to the continuous-time one. Fig.~\ref{fig5} gives the triggering time sequences of all agents.

\begin{figure}[!t]
  \centering
  \includegraphics[width=7.5cm]{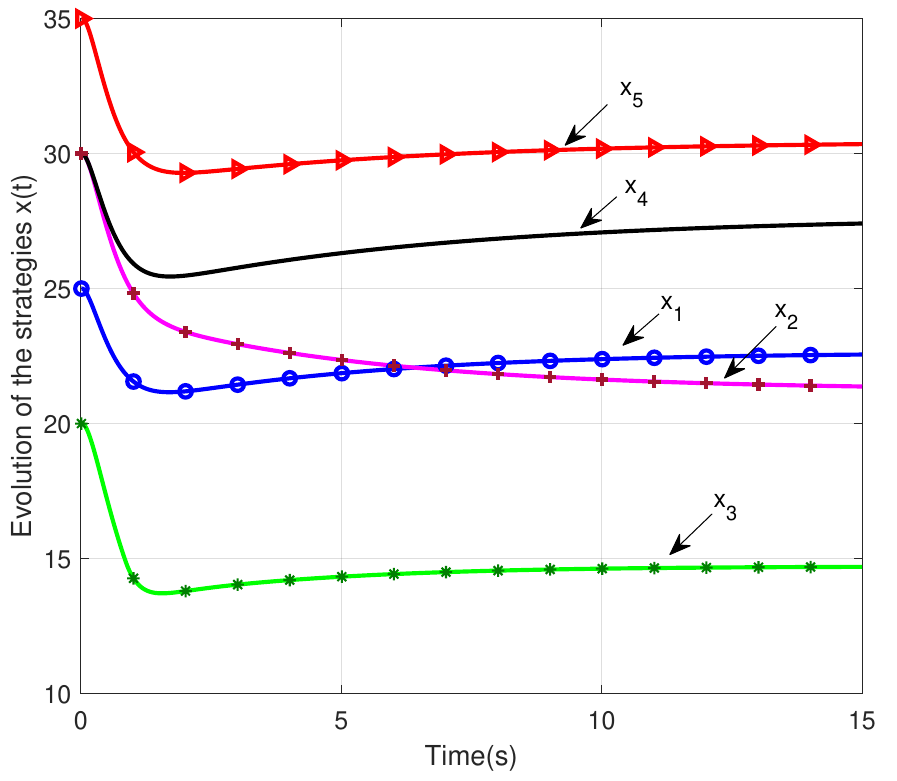}
  \caption{The evolution of the strategies of all agents following \eqref{dsur} with periodic communications}
  \label{fig3}
\end{figure}

\begin{figure}[!t]
  \centering
  \includegraphics[width=7.5cm]{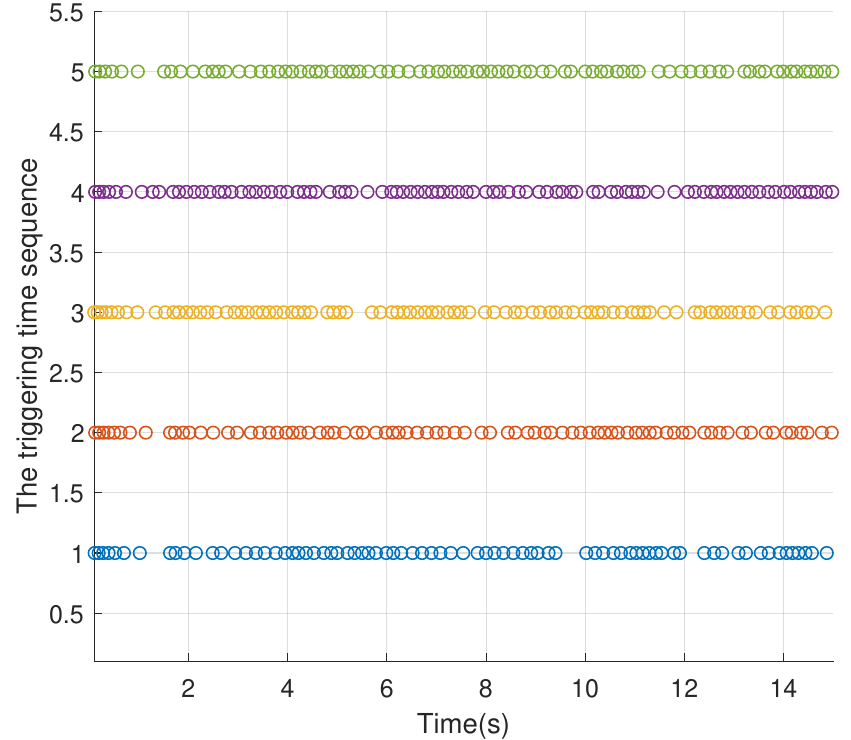}
  \caption{The triggering time sequences of the five agents}
  \label{fig5}
\end{figure}

In Fig.~\ref{fig3}, all agents' strategies, who follow the strategy-updating rule \eqref{dsur} in the periodic communication scheme with parameter $k=4$, $\alpha=5$, and $\triangle=0.1s$, can evolute to the NE of game $G$. From the comparison between Figs.~\ref{fig1} and \ref{fig3}, it is concluded that the event-triggered scheme results in a higher convergence rate. In addition, Table~\ref{tab2} shows the event times and event intervals of the five agents in the event-triggered scheme. It is seen that the the average event interval is greater than the sampling period and the communication frequency in the event-triggered system is less than that in the periodic one. According to the above simulation, the event-triggered communication scheme outperforms the periodic one in terms of convergence rates and communication frequencies.

\begin{table}[!t]
  \centering
  \caption{Event-triggered communication for the five agents}
  \label{tab2}
  \begin{tabular}{llllll}
    \hline
    Agent           & 1 & 2 & 3 & 4 & 5 \\
    \hline
    Event times         & 78 & 81 & 88 & 92 & 84 \\
    Min interval       & 0.05 & 0.05 & 0.04 & 0.05 & 0.05\\
    Mean interval     & 0.1906 & 0.1848 & 0.1686 & 0.1628 & 0.1783 \\
    Max interval      & 0.62 & 0.49 & 0.52 & 0.34 & 0.53\\
    \hline
  \end{tabular}
\end{table}

\section{Conclusions}
We have designed a distributed continuous-time strategy-updating rule for double-integrator agents, whose cost functions are continuous and not necessarily continuously differentiable. Our designed rule has been analyzed to ensure the evolution of agents' strategies to the NE of noncooperative games, if the communication graph is connected and undirected. This property is preserved in strongly connected and weight-balanced communication graphs. Then, discrete-time communication schemes for the implementation of the proposed rule are explored, such as periodic and event-triggered communication schemes. Furthermore, we have established the asymptotical convergence results in periodic and event-triggered communication schemes, and have taken care of the Zeno-behavior of the designed schemes. In future work, the influence of disturbances and time delays, and more complex agents' dynamics \cite{Li.2017,Sun.2017} can be considered in the model. And it may be an interesting issue to study the use of triggered communication schemes in the games with shared constraints.


%



\ifCLASSOPTIONcaptionsoff
  \newpage
\fi



%


\begin{thebibliography}{0}
\bibitem{Ratliff.2016}
 L. Ratliff, S. Burden, and S. Sastry, ``On the characterization of local Nash
equilibria in continuous games,'' \emph{IEEE Transactions on Automatic Control},
vol. 61, no. 8, pp. 2301-2307, Aug 2016.

\bibitem{Cortes.2015}
A. Cortes and S. Martinez, ``Slef-triggered bset-response dynamics for
continuous games,'' \emph{IEEE Transactions on Automatic Control}, vol. 60,
no. 4, pp. 1115-1120, Apr 2015.

\bibitem{Nisan.2007}
N. Nisan, T. Roughgarden, E. Tardos, and V. Vazirani, \emph{Algorithmic Game
Theory}.
England: Cambridge University Press, 2007.

\bibitem{Swenson.2018}
B. Swenson, R. Murry, and S. Kar, ``On best-response dynamics in
potential games,'' \emph{SIAM Journal on Control and Optimization}, vol. 56,
no. 4, pp. 2734-2767, 2018.

\bibitem{Chen.2019}
X. Chen, A. Brannstrom, and U. Dieckmann, ``Parent-preferred dispersal
promotes cooperation in structured populations,'' \emph{Prodeedings of the Royal
Society B-Biological Sciences}, vol. 286, no. 1895, pp. 1-8, Jan 2019.

\bibitem{Zheng.2018}
Y. Zheng, J. Ma, and L. Wang, ``Consensus of hybrid multi-agent systems,''
\emph{IEEE Transactions on Neural Networks and Learning Systems}, vol. 29,
no. 4, pp. 1359-1365, Apr 2018.

\bibitem{Qu.2018}
J. Qu, Z. Ji, C. Lin, and H. Yu, ``Fast consensus seeking on networks with
antagonistic interactions,'' \emph{Complexity}, pp. 1-15, 2018.

\bibitem{Yu.2018}
 J. Yu and Y. Shi, ``Scaled group consensus in multiagent systems
with first/second-order continuous dynamics,'' \emph{IEEE Transactions on
Cybernetics}, vol. 48, no. 8, pp. 2259-2271, Aug 2018.

\bibitem{Stankovic.2012}
M. Stankovic, K. Johansson, and D. Stipanovic, ``Distributed seeking
of Nash equilibria with applications to mobile sensor networks,'' \emph{IEEE
Transactions on Automatic Control}, vol. 57, no. 4, pp. 904-919, Apr 2012.

\bibitem{Ma.2019}
J. Ma, M. Ye, Y. Zheng, and Y. Zhu, ``Consensus analysis of hybrid
multiagent systems: a game-theoretic approach,'' \emph{International Journal of
Robust and Nonlinear Control}, vol. 29, no. 6, pp. 1840-1853, Apr 2019.

\bibitem{Gharesifard.2013}
B. Gharesifard and J. Cortes, ``Distributed convergence to Nash equilibria
in two-network zero-sum games,'' \emph{Automatica}, vol. 49, pp. 1683-1692,
Jun 2013.

\bibitem{Deng.2019}
 Z. Deng and S. Liang, ``Distributed algorithms for aggregative games of
multiple heterogeneous Euler-Lagrange systems,'' \emph{Automatica}, vol. 99, pp.
246-=0252, Jan 2019.

\bibitem{Ma.2013}
Z. Ma, D. Callaway, and I. Hiskens, ``Decentralized charging control
of large populations of plug-in electric vehicles,'' \emph{IEEE Transactions on
Control Systems Technology}, vol. 21, no. 1, pp. 67-78, Jan 2013.

\bibitem{Ibrahim.2018}
A. Ibrahim and T. Hayakawa, ``Nash equilibrium seeking with second-
order dynamics agents,'' in \emph{57th IEEE Conference on Decision and
Control}, Dec 2018, pp. 2514-2518.

\bibitem{Romano.2020}
A. R. Romano and L. Pavel, ``Dynamic NE seeking for multi-integrator
networked agents with disturbance rejection,'' \emph{IEEE Transactions on
Control of Network Systems}, vol. 7, no. 1, pp. 129-139, Mar 2020.

\bibitem{Bianchi.2019}
M. Bianchi and S. Grammatico, ``Continuous-time fully distributed
generalized Nash equilibrium seeking for multi-integrator agents,''
\emph{arXiv:1911.12266}, pp. 1-15, 2019.

\bibitem{Persis.2019c}
C. Persis and N. Monshizadeh, ``A feedback control algorithm to steer
networks to a Cournot-Nash equilibrium,'' \emph{IEEE Transactions on Control
of Network Systems}, vol. 6, no. 4, pp. 1486-1497, Dec 2019.

\bibitem{Zhang.2019}
Y. Zhang, S. Liang, X. Wang, and H. Ji, ``Distributed Nash equilibrium
seeking for aggregative games with nonlinear dynamics under external
disturbances,'' \emph{IEEE Transactions on Cybernetics}, vol. 50, no. 12, pp.
4876-4885, Dec 2019.

\bibitem{Huang.2020}
B. Huang, Y. Zou, and Z. Meng, ``Distributed-observer-based Nash
equilibrium seeking algorithm for quadratic games with nonlinear
dynamics,'' \emph{IEEE Transactions on Systems, Man, and Cybernetics:
Systems}, pp. 1-9, 2020.

\bibitem{Yuan.2018}
Y. Yuan, Z. Wang, and L. Guo, ``Event-triggered strategy design for
discrete-time nonlinear quadratic games with disturbance compensations:
the noncooperative case,'' \emph{IEEE Transactions on Systems, Man and
Cbernetics: Systems}, vol. 48, no. 11, pp. 1885-1896, Nov 2018.

\bibitem{Xue.2020}
S. Xue, B. Luo, and D. Liu, ``Event-triggered adaptive dynamic
programming for zero-sum game of partially unknown continuous-time
nonlinear systems,'' \emph{IEEE Transactions on Systems, Man, and Cybernetics:
Systems}, vol. 50, no. 9, pp. 3189-3199, Sep 2020.

\bibitem{Girard.2015}
A. Girard, ``Dynamic triggering mechanisms for event-triggered control,''
\emph{IEEE Transactions on Automatic Control}, vol. 60, no. 7, pp. 1992-1997,
Jul 2015.

\bibitem{Li.2016}
C. Li, X. Yu, W. Yu, T. Huang, and Z. Liu, ``Distributed event-triggered
scheme for economic dispatch in smart grids,'' \emph{IEEE Transactions on
Industrial Informatics}, vol. 12, no. 5, pp. 1775-1785, Oct 2016.

\bibitem{Liu.2019}
Q. Liu, M. Ye, J. Qin, and C. Yu, ``Event-triggered algorithms for
leaderfollower consensus of networked EulerLagrange agents,'' \emph{IEEE
Transactions on Systems, Man, and Cybernetics: Systems}, vol. 49, no. 7,
pp. 1435-1447, Jul 2019.

\bibitem{Sun.2020}
Y. Sun, Z. Ji, and K. Liu, ``Event-based consensus for general linear
multiagent systems under switching topologies,'' \emph{Complexity}, no. 5972749,
2020.

\bibitem{Ge.2020}
 X. Ge, Q. Han, L. Ding, Y. Wang, and X. Zhang, ``Dynamic event-triggered
distributed coordination control and its applications: A survey of trends
and techniques,'' \emph{IEEE Transactions on Systems, Man, and Cybernetics:
Systems}, vol. 50, no. 9, pp. 3112-3125, Sep 2020.

\bibitem{Ye.2017}
M. Ye and G. Hu, ``Distributed Nash equilibrium seeking by a consensus
based approach,'' \emph{IEEE Transactions on Automatic Control}, vol. 62, no. 9,
pp. 4811-4818, Sep 2017.

\bibitem{Gadjov.2019}
D. Gadjov and L. Pavel, ``A passivity-based approach to Nash equilibrium
seeking over networks,'' \emph{IEEE Transactions on Automatic Control}, vol. 64,
no. 3, pp. 1077-1092, Mar 2019.

\bibitem{Lu.2019}
 K. Lu, G. Jing, and L. Wang, ``Distributed algorithms for searching
generalized Nash equilibrium of noncooperative games,''
\emph{IEEE Transactions on Cybernetics}, vol. 49, no. 6, pp. 2362-2371, Jun 2019.

\bibitem{Godsil.2001}
 C. Godsil and G. Royle, \emph{Algebraic Graph Theory (Graduate Texts in
Mathematics)}.
New York, USA: Springer, 2001.

\bibitem{Dreves.2018}
 A. Dreves and M. Gerdts, ``A generalized Nash equilibrium approach
for optimal control problems of autonomous cars,'' \emph{Optimal Control
Applications and Methods}, vol. 39, pp. 326-342, Jan 2018.

\bibitem{Hobbs.2007}
B. Hobbs and J. Pang, ``Nash-Cournot equiliria in electric power markets
with piecewise linear demand functions and joint constraints,'' \emph{Operations
Research}, vol. 55, no. 1, pp. 113-127, Jan 2007.

\bibitem{Jafarpour.2011}
 S. Jafarpour, V. Cevher, and R. Schapire, ``A game theoretic approach
to expander-based compressive sensing,'' in \emph{2011 IEEE International
Symposium on Information Theory Proceedings}, 2011, pp. 464-468.

\bibitem{Yin.2011}
H. Yin, U. Shanbhag, and P. Mehta, ``Nash equilibrium problems with
scaled congestion costs and shared constraints,'' \emph{IEEE Transactions on
Automatic Control}, vol. 56, no. 7, pp. 1702-1708, Jul 2011.

\bibitem{Jean.2001}
 B. Jean and L. Claude, \emph{Fundamentals of Convex Analysis}.
Berlin, German: Springer, 2001.

\bibitem{Cortes.2008}
 J. Cortes, ``Discontinuous dynamical systems,'' \emph{IEEE Control Systems Magazine}, vol. 28, no. 3, pp. 36-73, Jun 2008.

\bibitem{Aubin.1984}
 J. Aubin, \emph{Differential Inclusions}.
Berlin: Springer, 1984.

\bibitem{Rosen.1965}
 J. Rosen, ``Existence and uniqueness of equilibrium points for concave
n-person games,'' \emph{Econometrica}, vol. 33, no. 3, pp. 520-534, Jul 1965.

\bibitem{Deng.2019b}
Z. Deng and X. Nian, ``Distributed generalized Nash equilibrium seeking
algorithm design for aggregative games over weight-balanced digraphs,''
\emph{IEEE Transactions on Neural Networks and Learning Systems}, vol. 30,
no. 3, pp. 695-706, Mar 2019.

\bibitem{Zeng.2019}
X. Zeng, J. Chen, S. Liang, and Y. Hong, ``Generalized Nash
equilibrium seeking strategy for distributed nonsmooth multi-cluster
game,'' \emph{Automatica}, vol. 103, pp. 20-26, May 2019.

\bibitem{Cai.2020}
X. Cai, F. Xiao, and B. Wei, ``A distributed strategy-updating rule
with event-triggered communication for noncooperative games,'' in \emph{39th
Chinese Control Conference}, Jul 2020, pp. 4747-4752.

\bibitem{Khalil.2002}
 H. Khalil, \emph{Nonlinear Systems}, 3rd ed.
USA: Prentice-Hall, 2002.

\bibitem{Basar.1999}
T. Basar and G. Olsder, \emph{Dynamic Noncooperative Game Theory}, 2nd ed.
USA: SIAM, 1999.

\bibitem{Persis.2019}
C. Persis and S. Grammatico, ``Distributed averaging integral Nash
equilibrium seeking on networks,'' \emph{Automatica}, vol. 110, no. 108548, Dec
2019.

\bibitem{Li.2017}
Y. Li, Y. Sun, and F. Meng, ``New criteria for exponential stability of
switched time-varying systems with delays and nonlinear disturbances,''
\emph{Nonlinear Analysis-Hybrid Systems}, vol. 26, pp. 284-291, Nov 2017.

\bibitem{Sun.2017}
 Y. Sun, Y. Tian, and X. Xie, ``Stabilization of positive switched linear
systems and its application in consensus of multiagent systems,'' \emph{IEEE
Transactions on Automatic Control}, vol. 62, no. 12, pp. 6608-6613, Aug
2017.
\end{thebibliography}

\end{document}